\documentclass[12pt,twoside]{article}

\usepackage{amsfonts}
\usepackage{amssymb}
\usepackage{amsmath}
\usepackage{graphicx}

\newcommand{\ignore}[1]{}

\def\@begintheorem#1#2{\par\bgroup{\sc #1\ #2. }\it\ignorespaces}
\def\@opargbegintheorem#1#2#3{\par\bgroup{\sc #1\ #2\ (#3). }\it\ignorespaces}
\def\@endtheorem{\egroup}
\newtheorem{theorem}{Theorem}[section]
\newtheorem{corollary}[theorem]{Corollary}
\newtheorem{lemma}[theorem]{Lemma}
\newtheorem{proposition}[theorem]{Proposition}
\newtheorem{example}[theorem]{Example}
\newtheorem{algorithm}[theorem]{Algorithm}
\newtheorem{definition}[theorem]{Definition}
\newcommand{\bt}[1]{\begin{theorem}\label{#1}}
\newcommand{\bc}[1]{\begin{corollary}\label{#1}}
\newcommand{\bl}[1]{\begin{lemma}\label{#1}}
\newcommand{\bp}[1]{\begin{proposition}\label{#1}}
\newcommand{\be}[1]{\begin{example}\rm\label{#1}}
\newcommand{\ba}[1]{\begin{algorithm}\rm\label{#1}}
\newcommand{\bd}[1]{\begin{definition}\rm\label{#1}}
\newcommand{\bpr}{\noindent {\em Proof. }}
\newcommand{\et}{\end{theorem}}
\newcommand{\ec}{\end{corollary}}
\newcommand{\el}{\end{lemma}}
\newcommand{\ep}{\end{proposition}}
\newcommand{\ee}{\end{example}}
\newcommand{\ea}{\end{algorithm}}
\newcommand{\ed}{\end{definition}}
\newcommand{\epr}{{\ \vbox{\hrule\hbox{%
\vrule height1.3ex\hskip0.8ex\vrule}\hrule}}\\\par}

\def\R{\mathbb{R}}
\def\Z{\mathbb{Z}}

\def \G {{\cal G}}
\def \l {\langle}
\def \r {\rangle}

\def \G {{\cal G}}

\def \l {\langle}
\def \r {\rangle}

\def \A {A^{(n)}}

\def \conv {{\rm conv}}

\def \supp {{\rm supp}}

\def \supp {{\rm supp}}

\def \zone {{\rm zone}}
\def \rank {{\rm rank}}

\begin{document}

\title{\bf Convex Integer Optimization
by Constantly Many Linear Counterparts}

\author{
Shmuel Onn
\thanks{Technion - Israel Institute of Technology, Haifa, Israel\newline
\texttt{onn@ie.technion.ac.il}}
\and
Michal Rozenblit
\thanks{Technion - Israel Institute of Technology, Haifa, Israel\newline
\texttt{michalro@tx.technion.ac.il}}
}

\date{\small Dedicated to the memory of Uri Rothblum}

\maketitle

\begin{abstract}
In this article we study convex integer maximization problems
with composite objective functions of the form $f(Wx)$, where
$f$ is a convex function on $\R^d$ and $W$ is a $d\times n$ matrix with small or
binary entries, over finite sets $S\subset \Z^n$ of integer
points presented by an oracle or by linear inequalities.

Continuing the line of research advanced by Uri Rothblum and his
colleagues on edge-directions, we introduce here the notion of
{\em edge complexity} of $S$, and use it to establish polynomial
and constant upper bounds on the number of vertices of the projection
$\conv(WS)$ and on the number of linear optimization
counterparts needed to solve the above convex problem.

Two typical consequences are the following.
First, for any $d$, there is a constant $m(d)$ such that the maximum number
of vertices of the projection of any matroid $S\subset\{0,1\}^n$
by any binary $d\times n$ matrix $W$ is $m(d)$ regardless of $n$ and $S$;
and the convex matroid problem reduces to $m(d)$ greedily solvable
linear counterparts. In particular, $m(2)=8$.
Second, for any $d,l,m$, there is a constant $t(d;l,m)$
such that the maximum number of vertices of the projection
of any three-index $l\times m\times n$ transportation polytope for any $n$
by any binary $d\times(l\times m\times n)$ matrix $W$ is $t(d;l,m)$;
and the convex three-index transportation problem
reduces to $t(d;l,m)$ linear counterparts solvable in polynomial time.
\end{abstract}

\newpage

\section{Introduction}

In this article we study convex integer maximization problems
and the closely related projections of the sets of feasible points.
Let $S\subset \Z^n$ be a finite set of integer points,
let $\conv(S)\subset\R^n$ be its convex hull, let $W$ be a
$d\times n$ integer matrix, and let $f:\R^d\rightarrow\R$
be a convex function. We study the problem of maximizing the
composite function $f(Wx)$ over $S$ and the projection of $\conv(S)$
by $W$ into $\R^d$, namely,
\begin{equation}\label{CDM}
\max\,\left\{f(Wx)\ :\ x\in S\right\}
\end{equation}
and
\begin{equation}\label{Projection}
\conv(WS)\ =\ \conv\{Wx\,:\,x\in S\}\ =\ \{Wx\,:\,x\in\conv(S)\}\ \subset\ \R^d\ .
\end{equation}

The sets $S$ we consider arise in two natural contexts.
First, in {\em combinatorial optimization}, in which case $S\subseteq\{0,1\}^n$
has some combinatorial structure and might be presented by
a suitable oracle. Second, in {\em integer programming}, where
\begin{equation}\label{IP}
S\ :=\ \left\{x\in\Z^n\,,\ Ax=b\,,\ l\leq x\leq u\right\}
\end{equation}
is the set of integer points satisfying a
given (standard) system of linear inequalities.

The optimization problem (\ref{CDM}) can also be interpreted as
a problem of multicriteria optimization, where each row of $W$ gives
a linear criterion $W_ix$ and $f$ compromises these criteria.
We therefore call $W$ the {\em criteria} matrix or {\em weight} matrix.

The projection polytope $\conv(WS)$ in (\ref{Projection}) and its vertices play a
central role in solving problem (\ref{CDM}): for any convex function $f$
there is an optimal solution $x$ whose projection $y:=Wx$ is a vertex of $\conv(WS)$.
In particular, the enumeration of all vertices of $\conv(WS)$ enables to compute
the optimal objective value for any given convex function $f$ by picking
that vertex attaining the best value $f(y)=f(Wx)$. So it suffices, and will be
assumed throughout, that $f$ is presented by a {\em comparison oracle}
that, queried on vectors $y,z\in\R^d$, asserts whether or not $f(y)<f(z)$.

This line of research was advanced by Uri Rothblum, to whom we
dedicate this article, and his colleagues, in several papers
including \cite{BHR,DHORW,HOR}, and culminated
in the edge-direction framework of \cite{OR}, see
also \cite[Chapter 2]{Onn}. In this article we continue this line of
investigation, and take a closer look on {\em coarse} criteria matrices;
that is, we assume that the entries of $W$ are small, presented
in unary, or even bounded by a constant and lie in $\{0,1,\dots,p\}$.
In multicriteria combinatorial optimization, this corresponds to the
weight $W_{i,j}$ attributed to element $j$ of the ground set $\{1,\dots,n\}$
under criterion $i$ being a small or even $\{0,1\}$ value for all $i,j$.

\vskip.2cm
Here is a typical result we obtain in convex combinatorial optimization,
where $S\subset\{0,1\}^n$ is the set of indicating vectors of bases of a
matroid over $\{1,\dots,n\}$.

\vskip.2cm\noindent{\bf Corollary \ref{Matroid_Plane}\ }
{\em The maximum number of vertices of $\conv(WS)$ for any $n$, any matroid
$S\subset\{0,1\}^n$ and any binary $2\times n$ matrix $W$, is $8$;
and $\max\{f(Wx):x\in S\}$ can be solved for any
convex $f$ by greedily solving $8$ linear counterparts over $S$.}

\vskip.2cm
Note that, in contrast, there are sets $S\subset\{0,1\}^n$ and $\{0,1\}$-valued
$2\times n$ matrices $W$ such that the projection $\conv(WS)$ has
$\Omega(\sqrt{n})$ vertices, see Theorem \ref{binary_sets_binary_matrices}.

More generally, we show in Section 4 that for any $d$ and $p$ there is a
constant $m(d,p)$ such that the maximum number of vertices of the projection
of any matroid $S\subset\{0,1\}^n$ by any $\{0,1,\dots,p\}$-valued $d\times n$
criteria matrix $W$ is $m(d,p)$ regardless of $n$ and $S$. Thus, $m(2,1)=8$
but the precise values of $m(d,p)$ for all higher $d$ and $p$ are unknown
and their determination poses a challenging open problem.

\vskip.2cm
Moving to convex integer programming, we have the following
typical result concerning multi-index transportation problems.

\vskip.2cm\noindent{\bf Corollary \ref{Multi-index}\ }
{\em For every fixed $d,p,l,m$, there exists a constant $t(d,p;l,m)$ such that,
for every $n$, every integer line-sums $a_{j,k},b_{i,k},c_{i,j}$,
and every $\{0,1,\dots, p\}$-valued $d\times(l\times m\times n)$ matrix $W$,
the projection $\conv(WS)$ of the set of $l\times m\times n$ tables
\begin{equation}\label{tables}
S\ :=\ \left\{x\in\Z_+^{l\times m\times n}\ :\ \sum_i x_{i,j,k}=a_{j,k}
\,,\ \sum_j x_{i,j,k}=b_{i,k}\,,\ \sum_k x_{i,j,k}=c_{i,j}\right\}\ ,
\end{equation}
has at most $t(d,p;l,m)$ vertices regardless of $n$ and the line sums.
Moreover, all these vertices can be enumerated, and for any convex $f$
the problem $\max\{f(Wx):x\in S\}$ solved, by solving $t(d,p;l,m)$
linear counterparts over $S$, and in polynomial time.}

\vskip.2cm
Other applications, to convex totally unimodular integer programs,
network flows and more generally transshipment problems with convex
multicriteria objective functions, and vector partition problems,
are discussed in Section 5.

These results and applications follow from a general
result which we now describe. Let $S\subset\Z^n$ be a finite set of integer
points. We define the {\em edge complexity} of $S$ to be the smallest
nonnegative integer $e(S)$ such that every edge of the polytope $\conv(S)$ is
parallel to some vector $v\in\Z^n$ with $\|v\|_1\leq e(S)$.
For a matrix $W$ we use the notation $\|W\|_\infty:=\max_{i,j} |W_{i,j}|$.
We establish the following theorem.

\vskip.2cm\noindent{\bf Theorem \ref{main}\ }
{\em Fix any $d$. Let $S\subset\Z^n$ be a finite set of integer points
and $W$ an integer $d\times n$ matrix. Then $\conv(WS)$ has
polynomially many $O\left((e(S)\cdot \|W\|_\infty)^{d(d-1)}\right)$ vertices.
Moreover, if $S$ is presented by a linear-optimization oracle, and
endowed with an upper bound $e$ on $e(S)$, then all vertices of $\conv(WS)$
can be enumerated, and $\max\{f(Wx)\,:\,x\in S\}$ solved for
every convex $f:\R^d\rightarrow\R$, in polynomial time.}

\vskip.2cm
The dependency in Theorem \ref{main} on
the unary size of $e(S)$ and $W$ cannot be relaxed: in
Theorems \ref{projection_of_binary_hulls} and \ref{projection_by_binary_matrices}
we show, respectively, that already projections into the plane $\R^2$ can have
exponentially many vertices, even when $S\subseteq\{0,1\}^n$ and hence the
edge complexity of $S$ satisfies $e(S)\leq n$, or when $W$ is a $\{0,1\}$-matrix.

Theorem \ref{main} continues the line of research taken in
\cite{BHR,BLMORWW,DHORW,HT,HOR,OR} implicitly or explicitly.
But in these papers, a complete set of edge-directions of $\conv(S)$ in $\R^n$
was required as part of the input, whereas here, only a bound on the edge
complexity is needed. Consequently we get constant bounds on the number of vertices and
number of linear counterparts needed, in Corollaries \ref{Matroid_Plane} and
\ref{Multi-index} and other applications, even in
situations where $\conv(S)$ has exponentially many edge-directions.

\vskip.2cm
The rest of the paper is organized as follows. In Section 2
we discuss edge complexity and projections of polytopes,
prove our general Theorem \ref{main} which leads to our polynomial
and constant bounds on the number of vertices of projections,
and describe the algorithm underlying it. In section 3 we construct
polytopes whose projections have large
number of vertices, in contrast with Theorem \ref{main}. In particular,
we show in Theorems \ref{projection_of_binary_hulls}
and \ref{projection_by_binary_matrices} respectively, that
projections into the plane of $\{0,1\}$-sets $S$ or by binary
matrices $W$ can have exponentially many vertices.

In Section 4 we focus on matroids in convex combinatorial optimization.
We show in Theorem \ref{Matroids_abridged} the existence of
bounds $m(d,p)$ on the number of vertices of the projection
by any $\{0,1,\dots,p\}$ valued $d\times n$ matrix $W$ of any
matroid $S$ over a ground set of any size $n$, and reduce the
convex matroid optimization problem to $m(d,p)$ linear counterparts.
We then restrict attention to binary matrices $W$ and give sharper
upper and lower bounds on $m(d):=m(d,1)$ in Theorem \ref{Bounds}.
For $d=2$ these coincide and we conclude Corollary \ref{Matroid_Plane}
with $m(2)=8$. We conclude Section 4 with a detailed discussion
of bicriteria matroid optimization in Example \ref{Example}.

In Section 5 we proceed to convex integer programming.
We discuss Graver bases and their connection to edge complexity.
In Section 5.1 we prove Theorem \ref{TUM} on
convex totally unimodular integer programs, and discuss applications
to convex multicriteria transshipment in Corollary \ref{transshipment}
and vector partition problems in Corollary \ref{partition}.
In Section 5.2 we prove Theorem \ref{N-fold} on convex $n$-fold
integer programming, discuss the application to muti-index
transportation problems and demonstrate Corollary \ref{Multi-index},
and prove an upper bound on the maximum number
$t(d;l,m):=t(d,1;l,m)$ of vertices in any projection of $3$-way
tables by binary criteria matrices in Proposition \ref{binary_table_upper_bound}.
We conclude in Section 6 with some open problems.

\vskip.2cm
Before proceeding we set some computational complexity
notation. The {\em binary length} of an integer $z$ is the number
$\l z\r=O(\log(|z|+1))$ of bits in its binary encoding. The binary length
$\l x\r$ of an integer vector $x$ is the sum of binary lengths
of its entries, and the length $\l E\r$ of a set $E$ of integer vectors
is the sum of lengths of its elements. Our algorithms have running times which are typically
polynomial in $n$ and $\|W\|_\infty$; in the binary length of the rest of
the data if any, which in integer programming is typically the binary length
$\l b,l,u \r$ of the right-hand side and lower and upper bounds; and they
perform a polynomial number of queries to the comparison oracle presenting
$f$ and the linear-optimization oracle presenting $S$ when relevant.

\section{Edge complexity and projections of polytopes}

We begin with some preliminaries on edge-directions and zonotopes.
A {\em direction} of an edge ($1$-dimensional face) $e$
of a polytope $P$ is any nonzero scalar multiple of $u-v$
with $u,v$ the vertices of $e$.
A {\em set of all edge-directions of $P$} is a set which contains
some direction of each edge of $P$.
The {\em normal cone} of a polytope $P\subset\R^n$ at its face $F$
is the relatively open cone of those linear functions $h\in\R^n$
maximized over $P$ precisely at points of $F$. A polytope $Z$ is a
{\em refinement} of a polytope $P$ if the normal cone of every vertex of $Z$
is contained in the normal cone of some vertex of $P$. The {\em zonotope}
generated by a set of vectors $E=\{e_1,\dots,e_m\}$ in $\R^d$ is the following
polytope, which is the projection by $E$ of the cube $[-1,1]^m$ into $\R^d$,
$$\zone(E)\ :=\
\conv\left\{\sum_{i=1}^m \lambda_i e_i\,:\,\lambda_i=\pm 1\right\}
\ =\ \sum_{i=1}^m \left[-e_i,e_i\right]\ \subset\ \R^d\ .$$
We have the following two lemmas, see for instance \cite{GS} or
\cite[Chapter 2]{Onn} and the references therein for more details.
The lemmas are illustrated in Figure \ref{refinement_figure} below.
\bl{Refinement} Let $P$ be a polytope and let $E$ be a finite
set of all edge-directions of $P$. Then the zonotope
$Z:=\zone(E)$ generated by $E$ is a refinement of $P$.
\el
\bl{Zonotope}
The number of vertices of any zonotope $Z:=\zone(E)$ generated by a set
$E$ of $m$ vectors in $\R^d$ is at most $2\sum_{k=0}^{d-1}{{m-1}\choose k}$.
For every fixed $d$, there is an algorithm that, given $E\subset\Z^d$,
outputs every vertex $u$ of $Z:=\zone(E)$ along with some $h_u\in\Z^d$
maximized over $Z$ uniquely at $u$, in time polynomial in $m$ and $\l E\r$.
\el
\begin{figure}[hbt]
\centerline{\includegraphics[scale=.48]{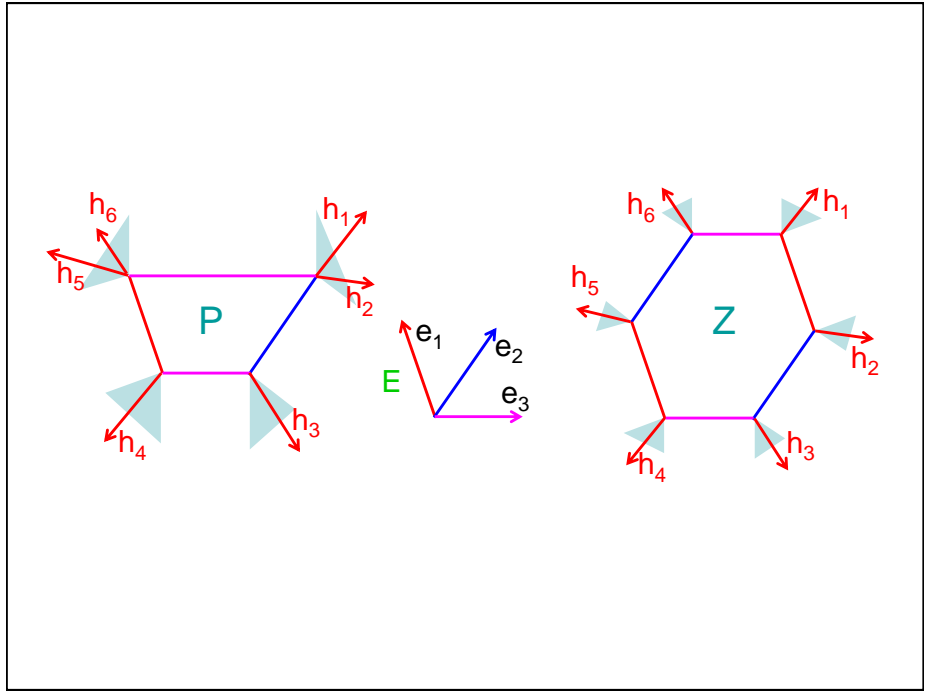}}
\caption{A zonotope refining a polytope and unique maximizers at their vertices }
\label{refinement_figure}
\end{figure}

We next introduce a new concept which is central to our development.
Let the {\em edge complexity} of a rational polytope $P$
be the smallest nonnegative integer $e(P)$ such that every
edge of $P$ is parallel to an integer vector $v$ with
$\|v\|_1\leq e(P)$. Let the edge complexity of a finite
set $S$ of rational points be $e(S):=e(\conv(S))$.
For matrix $W$ let $\|W\|_\infty:=\max_{i,j} |W_{i,j}|$.
A {\em linear-optimization oracle} for $S$
is one that, queried on $w\in\Z^n$, solves the
linear optimization problem $\max\{wx : x\in S\}$.

We are now ready to prove our theorem on edge complexity.
Its algorithmic statement is furnished by the
algorithm outlined following the proof of the theorem.

\bt{main}
Fix any $d$. Let $S\subset\Z^n$ be a finite set of integer points
and $W$ an integer $d\times n$ matrix. Then $\conv(WS)$ has
polynomially many $O\left((e(S)\cdot \|W\|_\infty)^{d(d-1)}\right)$ vertices.
Moreover, if $S$ is presented by a
linear-optimization oracle, and endowed with an upper bound $e$ on $e(S)$,
then all vertices of $\conv(WS)$ can be enumerated, and
$\max\{f(Wx)\,:\,x\in S\}$ solved for every convex $f:\R^d\rightarrow\R$,
in polynomial time.
\et
\bpr
Let $e:=e(S)$ and $q:=e\cdot \|W\|_\infty$. Let $P:=\conv(WS)$ be the projection
polytope. Now, consider any edge $f=[u,v]$ of $P$. Then there must be an edge $g$
of $\conv(S)$ whose projection $Wg$ is not a singleton and satisfies
$Wg\subseteq f$. Since $g$ is an edge of $\conv(S)$, it is parallel to
some $h\in\Z^n$ with $\|h\|_1\leq e(S)\leq e$. Therefore $Wh$ is an
edge-direction of $f$ and satisfies $Wh\in\{0,\pm 1,\dots,\pm q\}^d$,
since, for $i=1,\dots,d$,
$$|W_ih|\ \leq\ \sum_{j=1}^n|W_{i,j}||h_j|\ \leq\ \|W\|_\infty\sum_{j=1}^n|h_j|
\ =\ \|W\|_\infty \|h\|_1\ \leq\ q\ .$$
Since this is true for every edge of $P$, it follows that
$E:=\{0,\pm 1,\dots,\pm q\}^d$ is a set of all edge-directions of $P$.
Let $Z:=\zone(E)$ be the zonotope generated by the $(2q+1)^d$-set $E$.
Since $Z$ is a refinement of $P$ by Lemma \ref{Refinement}, the closure of
each normal cone of $P$ is the union of closures of normal cones of $Z$.
Therefore the number of vertices of $P$ is bounded by that of $Z$, which by
Lemma \ref{Zonotope} is at most
$$
2\sum_{k=0}^{d-1}{{(2q+1)^d-1}\choose k}
\ =\
2\sum_{k=0}^{d-1}{{(2e\cdot \|W\|_\infty+1)^d-1}\choose k}
\ =\
O\left((e\cdot \|W\|_\infty)^{d(d-1)}\right)\ .
$$
We continue with the algorithmic part of the theorem (see also
the description of the algorithm following this proof).
Let $e$ be the given upper bound, let again $q:=e\cdot \|W\|_\infty$
and $E:=\{0,\pm 1,\dots,\pm q\}^d$. By Lemma \ref{Zonotope},
in polynomial time we can produce every vertex $u$ of the zonotope
$Z:=\zone(E)$ along with some $h_u\in\Z^d$ maximized over $Z$ uniquely at $u$.
Now consider any vertex $v$ of $P:=\conv(WS)$.
Since $Z$ refines $P$ by Lemma \ref{Refinement}, there is a vertex
$u$ of $Z$ such that the normal cone of $Z$ at $u$ is contained
in the normal cone of $P$ at $v$. Therefore, $h_u$ is maximized over $P$
uniquely at $v$. Let $x_u\in S$ be a maximizer of $h_uW\in\Z^n$ over $S$.
Then, since for every $x\in\conv(S)$ and its image $y:=Wx\in\conv(WS)=P$
we have $h_uy=h_uWx$, it must be that $v=Wx_u$. So the algorithm
proceeds as follows (see description of the algorithm following this proof).
For each vertex $u$ of $Z$, use the linear-optimization oracle of $S$
to obtain a maximizer $x_u\in S$ of $h_uW\in\Z^n$ over $S$, and collect all such
$x_u$ in a set $X\subset\Z^n$ and all such $v=Wx_u$ in a set $V\subset\Z^d$.
Then $V$ is the set of vertices of $P$. Moreover, an optimal solution
$x^*$ is any point in $X$ attaining maximum $f(Wx)$, which
can be found using the comparison oracle of $f$.
\epr
The following furnishes the algorithm used and analyzed in the proof above.

\vskip.5cm
\noindent
\begin{tabular}{l}
\hline
\\
{\bf Projection and Convex Maximization Algorithm
(of  Theorem \ref{main})} \\
\\
{\bf Input:} Finite set $S\subseteq\Z^n$ presented by a linear-optimization oracle,\\
upper bound $e$ on its edge complexity $e(S)$, integer $d\times n$ matrix $W$,\\
and convex function $f:\R^d\rightarrow\R$ presented by a comparison oracle.\\
\\
1. {\bf Preprocessing:}\\
Let $q:=e\cdot \|W\|_\infty$. Construct the zonotope
$Z=\zone\{0,\pm 1,\dots,\pm q\}^d$, \\
and for each vertex $u$ of $Z$, some $h_u\in\Z^d$
maximized over $Z$ uniquely at $u$. \\
\\
2. {\bf Projection:}\\
Let $X:=\emptyset$, $V:=\emptyset$.
For each vertex $u$ of $Z$, obtain $x_u\in S$ maximizing \\
the linear function $h_uW$ over $S$ and
set $X:=X\cup\{x_u\}$, $V:=V\cup\{Wx_u\}$. \\
\\
3. {\bf Optimization:}\\
Obtain $x^*\in X$ attaining maximum value
of $f(Wx)$ over all $x\in X$.\\
\\
{\bf Output:}
The set $V$ of vertices of the projection $\conv(WS)$ \\
and the optimal solution $x^*$ of $\max\{f(Wx)\,:\,x\in S\}$. \\
\\
\hline
\end{tabular}

\section{Projections with exponentially many vertices }

In this section we construct sets $S$ of integer points and
weight matrices $W$ such that the projections $\conv(WS)$,
already into the plane $\R^2$, have many vertices,
contrasting our polynomial and constant bounds of the previous and next sections.

We begin by considering the case with both $S\subseteq\{0,1\}^n$
and $W\in\{0,1\}^{2\times n}$ binary valued. Then the projection satisfies
$WS\subseteq\{0,1,\dots,n\}^2$ and hence the number of vertices is $O(n^2)$ and
cannot be exponential. We show it can be $\Omega\left(\sqrt{n}\right)$.

\bt{binary_sets_binary_matrices}
For every $n$ there exist set $S\subseteq\{0,1\}^n$ and matrix
$W\in\{0,1\}^{2\times n}$ such that the projection of $\conv(S)$
by $W$ into the plane has more than $\sqrt{{1\over2}n}$ vertices.
\et
\bpr
Let $k:=\left\lfloor\sqrt{{1\over2}n}\right\rfloor$ so that $n\geq k+k^2$. Define a matrix
$W$ with each of its first $k$ columns being the unit vector ${\bf 1}_1\in\R^2$
and each of its next $k^2$ columns being the unit vector ${\bf 1}_2\in\R^2$.
Define a set $S$ in terms of the unit vectors ${\bf 1}_j\in\R^n$ by
$$S\ :=\ \left\{\sum_{j\leq i}{\bf 1}_j+
\sum_{k+1\leq j\leq k+i^2}{\bf 1}_j\ :\ 0\leq i\leq k\right\}\ .$$
Then $WS=\{(i,i^2)\,:\,0\leq i\leq k\}$ is a set of
$k+1>\sqrt{{1\over2}n}$ points in convex position, which comprise
the vertices of the projection of $\conv(S)$ by $W$ into the plane.
\epr
Next we construct several sets of integer points in polytopes whose
projections into the plane $\R^2$ have exponentially many vertices,
in contrast with Theorem \ref{main}.

First, we construct integer programs
defined by $\{-1,0,1\}$ matrices, and weight matrices
$W\in\Z^{2\times n}$ with binary size $\l W\r$ polynomial in $n$,
such that the integer hulls are $\{0,1\}$-polytopes, and
their projections have exponentially many vertices.
Note that the edge complexity of any $\{0,1\}$-polytope is at most $n$,
since any edge is parallel to the difference between its
two vertices which is a $\{-1,0,1\}$-valued vector.
\bt{projection_of_binary_hulls}
For every $k$ there are $m,n<2k^2$ and matrices $A\in\{0,\pm 1\}^{m\times n}$
and $W\in\Z_+^{2\times n}$ with $\log \|W\|_\infty <2k$, such that
the projection by $W$ of the $\{0,1\}$-polytope
$Q:=\conv\{x\in\Z^n\, :\, Ax=0\,,\, 0\leq x\leq 1\}$
into the plane $\R^2$ has exactly $2^k$ vertices.
\et
\bpr
Consider the following system of $m:=3{k\choose 2}$ equations
in $n:=k+4{k\choose 2}$ binary variables $x_i$, $0\leq i<k$ and
$x_{i,j},u_{i,j},v_{i,j},z_{i,j}$, $0\leq i<j<k$,
$$x_i-x_{i,j}-u_{i,j}=0\,,\quad x_j-x_{i,j}-v_{i,j}=0
\,,\quad x_i+x_j-x_{i,j}-z_{i,j}=0\,,\quad 0\leq i<j<k\ ,$$
and let $A$ be the corresponding $m\times n$ matrix.
Let $W$ be the matrix whose rows $w_1,w_2$ give the following
linear functions in the vector of variables $x:=(x_i,x_{i,j},u_{i,j},v_{i,j},z_{i,j})$,
$$w_1x\ :=\ \sum_{i=0}^{k-1} 2^i x_i\,,\quad
w_2x\ :=\ \sum_{i=0}^{k-1} 4^i x_i+2\sum_{0\leq i<j<k} 2^{i+j} x_{i,j}\ .$$
Clearly $m,n<2k^2$, $A$ is $\{0,\pm 1\}$-valued,
and $\l \|W\|_\infty \r<2k$, as claimed.

It is easy to see that the equations force that, for any $0\leq i<j<k$,
if $x_i=0$ or $x_j=0$ then $x_{i,j}=0$, whereas if $x_i=x_j=1$
then $x_{i,j}=1$. Thus, every choice of the partial variable
vector ${\bar x}:=(x_0,x_1,\dots,x_{k-1})\in\{0,1\}^k$ can be uniquely
extended to a vector $x\in\{0,1\}^n$ satisfying the above system of equations.
So the integer points in $Q$ are in bijection with such choices ${\bar x}$
and hence $Q$ has $2^k$ vertices. Moreover, for each such choice ${\bar x}$,
the equations force $x_{i,j}=x_ix_j$ for all $i,j$, and hence
\begin{eqnarray*}
(w_1x)^2 & = & \left(\sum_{i=0}^{k-1} 2^i x_i\right)^2\ =\
\sum_{i=0}^{k-1} (2^i)^2 x_i^2\, +\, 2\sum_{0\leq i<j<k} 2^i x_i 2^j x_j \\
& = & \sum_{i=0}^{k-1} 4^i x_i\, +\, 2\sum_{0\leq i<j<k}
2^{i+j}x_{i,j}\ =\ w_2x \ \ .
\end{eqnarray*}
When ${\bar x}$ runs over $\{0,1\}^k$, the value $w_1x$
runs over all values $0,1,\dots,2^k-1$. So the set
of integer points of $Q$ projects to the set
$\{(t,t^2)\ :\ t=0,1,\dots,2^k-1\}$ of $2^k$ points in convex
position, comprising the vertices of their convex hull $\conv(WS)$.
\epr
The last theorem of this section constructs integer programs defined
by $\{-1,0,1,2\}$-valued constraint matrices, and $\{0,1\}$-valued $2\times n$
weight matrices $W$, such that the number of vertices of the projections
of the integer hulls is exponential in $n$.
\bt{projection_by_binary_matrices}
For every $k$ there are $m,n<4k^4$,
matrices $A\in\{-1,0,1,2\}^{m\times n}$
and $W\in\{0,1\}^{2\times n}$, and lower and
upper bound vectors $l\in\{-\infty,0\}^n$,
$u\in\{1,\infty\}^n$, such that the integer hull
$Q:=\conv\{x\in\Z^n\, :\, Ax=0\,,\, l\leq x\leq u\}$
is a polytope whose projection by the binary matrix $W$
into the plane $\R^2$ has exactly $2^k$ vertices.
\et
\bpr
As in the proof of Theorem \ref{projection_of_binary_hulls},
consider again the following system of equations in binary
variables $x^0_i$, $0\leq i<k$ and
$x^0_{i,j},u_{i,j},v_{i,j},z_{i,j}$, $0\leq i<j<k$,
$$x^0_i-x^0_{i,j}-u_{i,j}=0\,,\quad x^0_j-x^0_{i,j}-v_{i,j}=0
\,,\quad x^0_i+x^0_j-x^0_{i,j}-z_{i,j}=0\,,\quad 0\leq i<j<k\ .$$
For each variable $x^0_i$ introduce $2i$ additional variables $x^r_i$
without bounds and $2i$ additional equations $2x^{r-1}_i-x^r_i=0$,
$r=1,\dots,2i$. Likewise, for each variable $x^0_{i,j}$ introduce $i+j+1$
additional variables $x^r_{i,j}$ without bounds and $i+j+1$
additional equations $2x^{r-1}_{i,j}-x^r_{i,j}=0$, $r=1,\dots,i+j+1$.
Then there are $n=$ variables and $m=$ equations with
$\{0,\pm 1,2\}$-valued defining matrix $A$.

Let $W$ be the $\{0,1\}$-valued matrix whose rows $w_1,w_2$
give the linear functions
$$w_1x\ :=\ \sum_{i=0}^{k-1} x^i_i\,,\quad
w_2x\ :=\ \sum_{i=0}^{k-1} x^{2i}_i+\sum_{0\leq i<j<k} x^{i+j+1}_{i,j}\ .$$

It is easy to see that the new equations force on the new variables the
relations $x^r_i=2^r x^0_i$ and $x^r_{i,j}=2^rx^0_{i,j}$ for all $i,j,r$.
Therefore, we obtain
$$w_1x\ =\ \sum_{i=0}^{k-1} 2^i x^0_i\,,\quad
w_2x\ :=\ \sum_{i=0}^{k-1} 4^i x^0_i+2\sum_{0\leq i<j<k} 2^{i+j} x^0_{i,j}\ .$$
The proof now proceeds as that of Theorem \ref{projection_of_binary_hulls},
that is, every choice of the partial variable
vector ${\bar x}:=(x^0_0,x^0_1,\dots,x^0_{k-1})\in\{0,1\}^k$
can be uniquely extended to a vector $x\in\Z^n$ satisfying
the system of equations, that is, an integer point in $Q$,
and these points project down onto the set
$\{(t,t^2)\ :\ t=0,1,\dots,2^k-1\}$ of $2^k$ points in convex
position, which comprise the vertices of their convex hull $\conv(WS)$.
\epr

\section{Convex matroid optimization}

In this section we study the convex multicriteria optimization
problem (\ref{CDM}) and the corresponding projection (\ref{Projection})
over matroids. We identify a matroid with its set of bases,
so we call $S\subset\{0,1\}^n$ a {\em matroid} if it is the set of
(indicating vectors of) bases of a matroid over $\{1,\dots,n\}$.
(We restrict attention to bases, but similar results hold for independent sets.)
In particular, our results below apply to spanning forests in graphs.
An {\em independence oracle} for $S$ is one that, queried on vector
$x\in\{0,1\}^n$, asserts whether or not $x$ is {\em independent}
in $S$, that is, whether or not $\supp(x)\subseteq \supp(z)$ for
some $z\in S$. As is well known, when a matroid $S$ is presented by an
independence oracle, the linear optimization problem $\max\{wx\,:\,x\in S\}$
can be easily solved by the greedy algorithm.
See \cite{Wel} for more details on matroids.

\vskip.2cm
As shown in Theorem \ref{binary_sets_binary_matrices}, when projecting
a binary set $S\subset\{0,1\}^n$ by a binary matrix $W\in\{0,1\}^{d\times n}$,
even to the plane $d=2$, the number of vertices can grow with $n$.
In contrast, we now show that for matroids,
the number of vertices is {\em constant}.

\bt{Matroids_abridged}
For every $d$ and $p$, the maximum number of vertices of $\conv(WS)$
for any $n$, any matroid $S\subset\{0,1\}^n$, and any $\{0,1,\dots,p\}$-valued
$d\times n$ matrix $W$, is a constant $m(d,p)$ which is independent of $n$,
$S$ and $W$. Moreover, if $S$ is presented by an independence oracle then
$\max\{f(WX):x\in S\}$ can be solved for any convex $f$ by greedily solving
$m(d,p)$ linear counterparts over $S$, and in polynomial time.
\et
\bpr
Let $S$ be any matroid. Then for every edge $e=[u,v]$ of the
matroid polytope $\conv(S)$ we have that $u-v={\bf 1}_i-{\bf 1}_j$ is equal
to the difference of some two unit vectors in $\R^n$; here is the proof,
included for completeness. Let $B_u:=\supp(u)$ and $B_v:=\supp(v)$ be the
corresponding bases. If $B_u\setminus B_v=\{i\}$ is a singleton then
$B_v\setminus B_u=\{j\}$ is a singleton as well in which case we are done.
Suppose then, indirectly, that this is not the case. Let $h\in\R^n$ be
uniquely maximized over $P$ at $e$, and pick an element $i$ in the symmetric
difference $B_u\Delta B_v:=(B_u\setminus B_v)\cup(B_v\setminus B_u)$ of
minimum value $h_i$. Without loss of generality assume $i\in B_u\setminus B_v$.
Then there is a $j\in B_v\setminus B_u$ such that $B:=B_u\setminus\{i\}\cup\{j\}$
is also a basis. Let $x:={\bf 1}_B\in S$ be the indicator of $B$.
Now $|B_u\Delta B_v|>2$ implies that $x\neq u$ and $x\neq v$.
By the choice of $i$ we have $hx=hu-h_i+h_j\geq hu$. So $x$ is also a
maximizer of $h$ over $P$ and hence $x\in e$. But no $\{0,1\}$-vector
is a convex combination of others, a contradiction.

Thus, every edge of $\conv(S)$ is parallel to some ${\bf 1}_i-{\bf 1}_j$
and hence $S$ has constant edge complexity $e(S)=2$. Also, the criteria matrix
satisfies $\|W\|_\infty=p$. Finally, a linear-optimization oracle for $S$ is realized
from the independence oracle presenting it by the greedy algorithm.
The theorem now follows from Theorem \ref{main}.
\epr

As mentioned in the introduction, the determination of the constants
$m(d,p)$ is a challenging open problem. We now restrict attention to
{\em binary} criteria matrices, that is, $\{0,1\}$-valued $d\times n$ matrices
$W$. Let $m(d):=m(d,1)$ be the corresponding maximum number of vertices
of the projection $\conv(WS)$ of any matroid $S$ by any binary matrix.
Introduce the zonotope $Z(d):=\zone\{-1,0,1\}^d$ and let $z(d)$ be the number
of vertices of $Z(d)$. We next provide an upper bound on $z(d)$ and $m(d)$.

\bl{binary_upper_bound}
The maximum number of vertices of the projection $\conv(WS)$ into $\R^d$ of any
matroid $S$ by any binary $d\times n$ matrix $W$ satisfies the upper bound
$$m(d)\ \leq\ z(d)\ \leq\
2\sum_{k=0}^{d-1}{{1\over2}(3^d-3)\choose k}\ .$$
\el
\bpr
The proof of Theorem \ref{main} implies that every edge of
$\conv(WS)$ is parallel to $Wh$ for some edge-direction $h$ of $\conv(S)$.
Since $S$ is a matroid, the proof of Theorem \ref{Matroids_abridged} implies
that every edge of $\conv(S)$ is parallel to some ${\bf 1}_i-{\bf 1}_j$.
Since $W\in\{0,1\}^{d\times n}$, it follows that a set of
all edge-directions of $\conv(WS)$ is given by
$$\{Wh\,:\,h={\bf 1}_i-{\bf 1}_j\,,\ 1\leq i<j\leq n\}
\ \subseteq\ \{-1,0,1\}^d\ .$$
So $Z(d)=\zone\{-1,0,1\}^d$
refines $\conv(WS)$ and the latter has at most $m(d)\leq z(d)$
vertices. Now, $Z(d)$ is homothetic to  the zonotope generated by the
following set containing
only one of each antipodal pair of nonzero vectors in $\{-1,0,1\}^d$,
$$\{v\in\{-1,0,1\}^d\ :\ v\neq 0\
\mbox{and the first nonzero entry of $v$ is $1$}\}\ .$$
Since this set has ${1\over2}(3^d-1)$ elements, the bound now
follows from Lemma \ref{Zonotope}.
\epr
We proceed to develop some lower bounds on $m(d)$ (and hence also
on $z(d)$). The {\em uniform matroid} of rank $r$ and order $n$
is the matroid having all $r$-subsets of $\{1,\dots,n\}$ as bases.
We identify it with the set $S^n_r\subset\{0,1\}^n$
of all vectors $x\in\{0,1\}^n$ with $|\supp(x)|=r$.
Let $u(d)$ be the maximum number of vertices of the projection
$\conv(WS^n_r)$ of any uniform matroid under any binary $d\times n$ matrix $W$.
Clearly, $u(d)\leq m(d)$ for all $d$. For positive integers $k,d$,
let $W_d^k$ be the $d\times k2^d$ binary criteria matrix whose $k2^d$ columns
consist of $k$ copies of each vector in $\{0,1\}^d$.
Finally, let $u_r^k(d)$ be the number of vertices of
$\conv(W_d^kS^{k2^d}_r)$; note that it is well defined and independent
of the order of columns of $W_d^k$ by the symmetry of uniform matroids.
Clearly, $u_r^k(d)\leq u(d)\leq m(d)$ for all $k,r,d$.
Beyond providing lower bounds on $m(d)$, the determination of
$u(d)$ and $u_r^k(d)$ is an interesting problem in its own right;
the following three statements characterize the vertices of
$\conv(W_d^kS^{k2^d}_r)$ and determine $u_r^k(d)$, respectively,
for $k\geq r$, $k=1$ and $r=2$, and $k\geq 2$ and $r=k+1$, for all $d$.

\bp{k=r}
For $k\geq r$, $\conv(W^k_d S^{k2^d}_r)=[0,r]^d$ is a cube
and hence $u_r^k(d)=2^d$.
\ep
\bpr
Clearly the image of $S^{k2^d}_r$ under $W^k_d$ satisfies
$W^k_d S^{k2^d}_r\subseteq\{0,1,\dots,r\}^d$ and hence is contained in the cube.
On the other hand, every vertex $v$ of the cube $[0,r]^d$ is a multiple
$v=r\cdot u$ for some $u\in\{0,1\}^n$; letting $J$ be any set of $r\leq k$
distinct indices of columns of $W^k_d$ which are equal to $u$, and
letting $x:=\sum_{j\in J}{\bf 1}_j\in S^{k2^d}_r$, we find that
$v=W^k_d x\in W^k_d S^{k2^d}_r$ is in the image. So the image is contained
in $[0,r]^d$ and contains all its vertices and hence
$\conv(W^k_d S^{k2^d}_r)=[0,r]^d$ as claimed.
\epr

We assume next $k<r$ and begin with $k=1$ and $r=2$.
Let $V:=\{0,1,2\}^d$ and for $i=0,1,\dots,d$ let $V_i$ be the set of
vectors in $V$ with precisely $i$ entries equal to $1$.

\bp{k=1_r=2}
For $k=1$ and $r=2$, the image satisfies $W^1_d S^{2^d}_2=V\setminus V_0$,
the set of vertices of its convex hull $\conv(W^1_d S^{2^d}_2)$ is equal to
$V_1$, and therefore $u_2^1(d)=d2^{d-1}$.
\ep
\bpr
Let $n:=2^d$. The image $W^1_d S^n_2$ consists of all sums $u+v$ of distinct
pairs $u,v\in\{0,1\}^d$. Each such sum $u+v$ is in $V$ and
for some $i$ we have $u_i\neq v_i$ so $u_i+v_i=1$ and hence
$u+v\not\in V_0$. Conversely, consider $z\in V\setminus V_0$ with
$z=\sum_{i\in I}{\bf 1}_i+\sum_{j\in J}2{\bf 1}_j$.
Then $z=u+v$ with $u:=\sum_{i\in I\cup J}{\bf 1}_i$
and $v:=\sum_{j\in J}{\bf 1}_j$ distinct since $I\neq\emptyset$.

We proceed with the claim about vertices. Consider any
$v\in W^1_d S^n_2\setminus V_1$. By what we just proved,
$v\not\in V_0$ and so $v_i=1$ for some $i$. Since $v\not\in V_1$
we have that $v+{\bf 1}_i$ and $v-{\bf 1}_i$ are not in $V_0$ hence
in the image $W^1_d S^n_2$, and
$v={1\over 2}\left((v+{\bf 1}_i)+(v-{\bf 1}_i)\right)$ is a convex combination
of them and hence not a vertex. Conversely, consider any $v\in V_1$. Then $v_i=1$
in precisely one entry. Define a vector $h\in\Z^d$ by
$$h_j\ :=\ \left\{\begin{array}{rl}
      -1, & \text{if } v_j=0,\\
       0, & \text{if } v_j=1,\\
       1, & \text{if } v_j=2.
    \end{array}\right.$$
Then $hu\leq hv$ for all $u\in V$, with equality
if and only if $u\in\left\{v-{\bf 1}_i\,,v\,,v+{\bf 1}_i\right\}$.
But $v-{\bf 1}_i\,,v+{\bf 1}_i\in V_0$ so $v$ is the unique
maximizer of $h$ over $W^1_d S^n_2$ hence a vertex.

So the number of vertices of $\conv(W^1_d S^n_2)$ is
$u_r^k(d)=|V_1|=d2^{d-1}$ as claimed.
\epr

We proceed with the case of $r=k+1$ for any $k\geq 2$.
Let $V:=\{0,1,\dots,r\}^d$ and for $i=0,1,\dots,d$ let $V_i$
be the set of vectors in $V$ with $i$ entries not equal to $0$ of $r$.

\bl{r=k+1}
For $k\geq 2$ and $r=k+1$, the image satisfies $W^k_d S^{k2^d}_r=V\setminus V_0$,
the set of vertices of $\conv(W^k_d S^{k2^d}_r)$ consists of those vectors
$v\in V_1$ which have one entry equal to $1$ or $r-1$ and all other entries
equal to $0$ or $r$, and therefore $u_r^k(d)=d2^d$.
\el
\bpr
Let $n:=k2^d$. The image $W^k_d S^n_r$ consists of all sums $z=\sum_{i=1}^r u^i$
of $r$ vectors $u^i\in\{0,1\}^d$ with at least two of them distinct.
Clearly each such sum $z$ is in $V$ and for some $i,j,h$ we have
$u^i_h\neq u^j_h$ so $u^i_h+u^j_h=1$ and hence $1\leq z_h\leq r-1$
so $z\not\in V_0$. Conversely, consider $z\in V\setminus V_0$.
Then $z=\sum_{i=1}^r\sum_{j\in I_i}i{\bf 1}_j$ for some pairwise
disjoint $I_i$ with $I_h\neq\emptyset$ for some $1\leq h\leq r-1$.
For $i=1,\dots,r$ let $u^i:=\sum_{j\in I_i\cup\cdots\cup I_r}{\bf 1}_j$.
Then $u^h\neq u^{h+1}$ and $z=\sum_{i=1}^r u^i$ so $z\in W^k_d S^n_r$.

We proceed with the claim about vertices. First we claim that the set of
vertices of $v\in W^k_d S^n_r$ is contained in $V_1$. Consider any
$v\in W^k_d S^n_r\setminus V_1$. By what we just proved, $v\not\in V_0$ and
so $1\leq v_i\leq r-1$ for some $i$. Since $v\not\in V_1$ we have that
$v+{\bf 1}_i$ and $v-{\bf 1}_i$ are not in $V_0$ hence in the image
$W^k_d S^n_r$, and $v={1\over 2}\left((v+{\bf 1}_i)+(v-{\bf 1}_i)\right)$
is a convex combination of them and hence not a vertex.

Next, consider any $v\in V_1$. Let $i$ be the unique index where
$1\leq v_i\leq r-1$, and for $t=0,1,\dots,r$ let $v^t$ be the vector defined by
$$v^t_j\ :=\ \left\{\begin{array}{rl}
      v_j, & \text{if } j\neq i,\\
       t, & \text{if } j=i.
    \end{array}\right.$$
Define a vector $h\in\Z^d$ by
$$h_j\ :=\ \left\{\begin{array}{rl}
      -1, & \text{if } v_j=0,\\
       0, & \text{if } v_j=1,\\
       1, & \text{if } v_j=2.
    \end{array}\right.$$
Then $hu\leq hv$ for all $u\in V$, with equality
if and only if $u\in U:=\left\{v^0,v^1,\dots,v^r\right\}$.
Therefore $W^k_d S^n_r\cap U=\left\{v^1,\dots,v^{r-1}\right\}$
spans a $1$-dimensional face of $\conv(W^k_d S^n_r)$ with vertices
$v^1,v^{r-1}$, which are therefore also vertices $\conv(W^k_d S^n_r)$.
Thus, the vertices consists precisely of all vectors $v\in V_1$ having
one entry $1$ or $r-1$. Therefore the number of vertices
of $\conv(W^k_d S^n_r)$ is $u_r^k(d)=d2^d$ as claimed.
\epr
Combining Lemma \ref{binary_upper_bound} and Lemma \ref{r=k+1}
we obtain at once the following theorem.
\bt{Bounds}
The maximum number of vertices of the projection polytope $\conv(WS)$ of any matroid
$S\subset\{0,1\}^n$ by any $\{0,1\}$-valued $d\times n$ criteria matrix $W$ satisfies
$$d2^d\ \leq\ m(d)\ \leq\
2\sum_{k=0}^{d-1}{{1\over2}(3^d-3)\choose k}\ =\ O\left(3^{d(d-1)}\right)\ .$$
For $d\leq 2$ the lower and upper bounds coincide
and we have $m(1)=2$ and $m(2)=8$.
\et

Note that already for $d=3$ and $d=4$ the bounds are far apart,
$24\leq m(3)\leq 158$ and $64\leq m(4)\leq 19840$, and we do not know the
exact values. The value $m(2)=8$ and Theorem \ref{Matroids_abridged}
give the following corollary mentioned in the introduction.

\bc{Matroid_Plane}
The maximum number of vertices of $\conv(WS)$ for any $n$, any matroid
$S\subset\{0,1\}^n$ and any binary $2\times n$ matrix $W$, is $8$;
and $\max\{f(Wx):x\in S\}$ can be solved for any
convex $f$ by greedily solving $8$ linear counterparts over $S$.
\ec

\begin{figure}[hbt]
\centerline{\includegraphics[scale=.48]{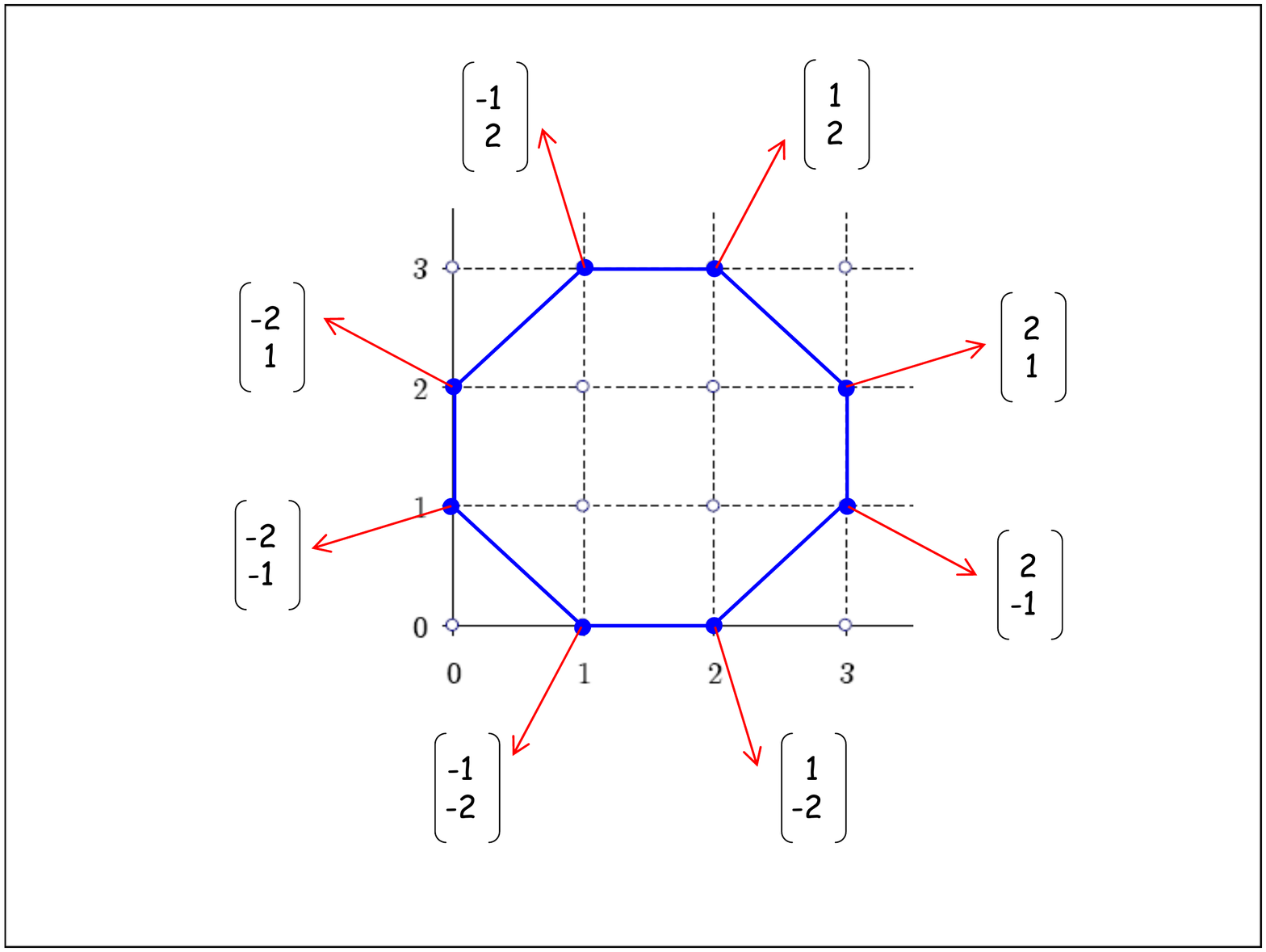}}
\caption{The octagon $\conv(W^2_2S^8_3)$ and homothety of zonotope $Z(2)$}
\label{octagon}
\end{figure}

We conclude this section with an explicit description of the matroid
and zonotope realizing the matching lower and upper bounds on $m(2)=8$,
and the linear counterparts used in solving any
convex bicriteria matroid optimization problem.

\be{Example}
Let $d=2$, $k=2$, $r=3$, and $n=8$, let $S^8_3$ be the uniform matroid\break
$S^8_3=\{(1,1,1,0,0,0,0,0),\dots,(0,0,0,0,0,1,1,1)\}$,
and let $W^2_2$ be the matrix
$$W^2_2\ :=\
\left(
\begin{array}{cccccccc}
  0 & 0 & 1 & 1 & 0 & 0 & 1 & 1 \\
  0 & 0 & 0 & 0 & 1 & 1 & 1 & 1 \\
\end{array}
\right)\,.
$$
Then, as predicted by Lemma \ref{r=k+1}, $\conv(W^2_2S^8_3)$ is an octagon,
with vertex set
$$U\ :=\ \{(1,0),(2,0),(3,1),(3,2),(2,3),(1,3),(0,2),(0,1)\}\ ,$$
homothetic to the zonotope
$Z(2)=\zone\{-1,0,1\}^2$, confirming $m(2)=z(2)=8$.
A set of $8$ vectors $H$ which contains, for every vertex $u\in U$,
a linear function $h_u$ which is maximized over $\conv(W^2_2S^8_3)$
and $Z(2)$ uniquely at $u$, is provided by
$$H=\{(-1,-2),(1,-2),(2,-1),(2,1),(1,2),(-1,2),(-2,1),(-2,-1)\}\ .$$
The polytope $\conv(W^2_2S^8_3)$ and zonotope $Z(2)$, along with the
set $H$, are depicted in Figure \ref{octagon}. Given now {\em any} $n$,
{\em any} matroid $S\subset\{0,1\}^n$ presented by an independence oracle,
{\em any} binary $2\times n$ matrix $W,$ and {\em any} convex function
$f:\R^2\rightarrow\R$ presented by a comparison oracle,
to solve $\max\{f(Wx)\,:\,x\in S\}$ using the algorithm
of Theorem \ref{main},
proceed as follows: for each of the $8$ vectors $h\in H$ use the greedy
algorithm to find $x_h\in S$ maximizing $hW$ over $S$; output that $x_h$
attaining best value $f(Wx_h)$.
\ee

\section{Convex integer programming}

We now turn to convex integer programming, where the set of feasible points is
\begin{equation}\label{IP_set}
S\ :=\ \left\{x\in\Z^n\,,\ Ax=b\,,\ l\leq x\leq u\right\}.
\end{equation}
We need to define the Graver basis of an integer matrix,
see \cite{Onn} for more details.
Define a partial order $\sqsubseteq$ on $\R^n$ where $x\sqsubseteq y$
if $x_iy_i\geq 0$ (that is, $x$ and $y$ lie in the same orthant) and
$|x_i|\leq |y_i|$ for $i=1,\ldots,n$. The {\em Graver basis} of an
integer matrix $A$ is defined to be the finite set $\G(A)\subset\Z^n$
of $\sqsubseteq$-minimal elements in $\{x\in\Z^n: Ax=0,\ x\neq 0\}$.
For example, the Graver basis of $A:=(1\,\ 2\,\ 1)$ is the set
$$\G(A)\ =\ \pm\left\{(2,-1,0),(0,-1,2),(1,0,-1),(1,-1,1)\right\}\ .$$
Let
$$\|\G(A)\|_1\ :=\ \max\left\{\|x\|_1\,:\,x\in\G(A)\right\}\ .$$

\vskip.2cm
We have the following Lemma.
\bl{Graver_projections}
Let $S$ be as in (\ref{IP_set}) above. Then $\|\G(A)\|_1$ is an upper bound on $e(S)$.
\el
\bpr
Let $S$ be as in (\ref{IP_set}). It has been shown in \cite{DHORW} that the
Graver basis of $A$ is a set of all edge-directions of $\conv(S)$;
here is the proof, included for completeness.
Consider any edge $e=[x,z]$ of $\conv(S)$ and
let $h:=z-x$. Then $Ah=0$ and hence $h=\sum_i h_i$ is a sum of some elements
$h_i\in\G(A)$, all lying in the same orthant, see \cite[Chapter 3]{Onn}. We claim
that $x+h_i\in \conv(S)$ for all $i$. Indeed, $h_i\in\G(A)$ implies  $A(x+h_i)=b$,
and $l\leq x,x+h \leq u$ and $h_i\sqsubseteq h$ imply $l\leq x+h_i\leq u$.

Now let $w\in\Z^n$ be uniquely maximized over $\conv(S)$ at
the edge $e$. Then $wh_i=w(x+h_i)-wx\leq 0$ for all $i$.
But $\sum wh_i=wh=wz-wx=0$, implying that in fact $wh_i=0$
and hence $x+h_i\in e$ for all $i$. This implies that $h_i$
is a direction of $e$ (in fact, all $h_i$ are the same
and $h$ is a multiple of some Graver basis element).

So any edge of $\conv(S)$ is parallel to some $v\in\G(A)$
and hence $e(S)\leq \|\G(A)\|_1$.
\epr

In the next two subsections we discuss two important classes of integer programs and
their applications, where the multicriteria convex problem reduces to a small or
constant number of linear counterparts, and can be solved altogether in polynomial
time, using Theorem \ref{main} together with Lemma \ref{Graver_projections}
and some more ingredients.

\subsection{Convex totally unimodular integer programs}

We have the following theorem on convex
totally unimodular integer programming.

\bt{TUM}
Fix any $d$. Consider the convex integer programming problem
\begin{equation}\label{TUM_problem}
\max\,\left\{f(Wx)\ :\ x\in S \right\}\,,\quad
S:=\left\{x\in\Z^n\,,\ Ax=b\,,\ l\leq x\leq u\right\},
\end{equation}
with $A$ a totally unimodular $m\times n$ matrix,
$b\in\Z^m$, $l,u\in\Z^n$, $W$ an integer $d\times n$ matrix,
and $f:\R^d\rightarrow\R$ a convex function. Then $\conv(WS)$ has
polynomially many $O\left((n\cdot \|W\|_\infty)^{d(d-1)}\right)$ vertices.
Further, all vertices of $\conv(WS)$ can be enumerated and
the convex integer program solved in time polynomial in $\|W\|_\infty$
and $\l b,l,u\r$.
\et
\bpr
A {\em circuit} of an integer matrix $A$ is a primitive linear dependency
on the columns of $A$, that is, a nonzero integer vector satisfying
$Ax=0$ whose support $\supp(x)$ is minimal under inclusion and whose
entries are relatively prime integers. It then follows from Cramer's rule
that every circuit $x$ of $A$ satisfies $\|x\|_1\leq(r+1)\Delta(A)$, where
$r=\rank(A)\leq \min\{m,n\}$ is the rank of $A$ and $\Delta(A)$ denotes
the maximum absolute value of a determinant of a square submatrix of $A$.

Now, for a totally unimodular $A$, the Graver basis $\G(A)$ consists
precisely of the circuits of $A$, see \cite[Lemma 3.19]{Onn}, and $\Delta(A)=1$.
Therefore for totally unimodular $A$ we obtain $\|\G(A)\|_1\leq n+1$ and hence
by Lemma \ref{Graver_projections} also $e(S)\leq n+1$.
Moreover, if $A$ is totally unimodular then we can realize in
polynomial time a linear optimization oracle for $S$ since we can
optimize linear functions over $S$ in polynomial time using linear programming.
The theorem therefore now follows from Theorem \ref{main}.
\epr
A specially appealing application of Theorem \ref{TUM} is to convex
multicriteria integer transshipment and in particular network flow problems.
In the integer transshipment problem, the feasible set is $S$ as in
(\ref{TUM_problem}), where $A$ is the totally unimodular
vertex-arc incidence matrix of the directed graph $G=(V,E)$ underlying the
problem; $b\in\Z^V$ is the demand vector, with $v$ a supplier if $b_v$
is negative and a consumer if $b_v$ is positive; $l_e$ and $u_e$ are lower
and upper capacities on flow on arc $e$; $W$ is the criteria matrix with $W_ix$
the value of transshipment $x\in\Z^E$ under criterion $i$; and $f$ compromises
these criteria. We have the following immediate corollary of Theorem \ref{TUM}.

\bc{transshipment}
Fix any $d$. Then for every directed graph $G$, demand $b$, lower and
upper capacities $l,u$, criteria matrix $W$, and convex function $f$
presented by a comparison oracle, the convex $d$-criteria integer transshipment
problem (\ref{TUM_problem}) can be solved in time polynomial in
$\|W\|_\infty$ and the binary size $\l b,l,u\r$ of the rest of the data.
\ec

Theorem \ref{TUM} can be obtained also through the results of \cite{BLOW}.
However, therein, one needs to repeatedly solve
extended linear programs of the form
$$\{x\in\R^n\,,\ Ax=b\,,\ Wx=u\,,\ l\leq x\leq u\}$$
over the so called {\em fibers} of various $u\in\Z^d$; so the
totally unimodular or specific combinatorial structure of the
matrix $A$ is lost. For instance, in solving the convex transshipment problem
in Corollary \ref{transshipment} using our method here, the linear
optimization oracle is very efficiently realizable by suitable
network flow algorithms, whereas the derivation through \cite{BLOW}
requires the use of a generic linear programming solver.

We conclude this subsection with an application to the vector partition problem
studied in \cite{BHR,HOR}. In this problem we need to partition $n$ given items
among $p$ players so that player $i$ gets $b_i$ items. We are given
$p\times n$ integer utility matrices $W^1,\dots,W^d$, where $W^k_{i,j}$ is the utility of
item $j$ to player $i$ under criterion $k$. A partition is naturally represented
by a $p\times n$ matrix $X$ with $X_{i,j}=1$ if item $j$ is assigned
to player $i$ and $X_{i,j}=0$ otherwise.
Let $W^k\cdot X:=\sum_{i,j}W^k_{i,j}\cdot X_{i,j}$ be the total utility under
criterion $k$ of partition $X$. The goal is to find
a partition $X$ attaining maximum value $f(W^1\cdot X,\dots,W^d\cdot X)$
where $f:\R^d\rightarrow\R$ is a convex function.

Since the partitions form the feasible points of the totally unimodular system
$$S\ =\ \{X\in\Z_+^{p\times n}\,:\,\ \sum_i X_{i,j}=1\,,\ \sum_j X_{i,j}=b_i\}\ ,$$
we obtain the following immediate corollary of Theorem \ref{TUM}.

\bc{partition}
Fix any $d$. Then the partition problem with any number $n$ of items,
number $p$ of players, and $b\in\Z^p$ with $p\leq n= \sum b_i$, and
$p\times n$ matrices $W^1,\dots,W^d$ and convex function $f$,
can be solved in time polynomial in $n$ and
$\|W^1\|_\infty,\dots,\|W^d\|_\infty$.
\ec

We emphasize that here the number $p$ of players can be {\em variable}, whereas in
previous work (see \cite{HOR} and references therein) it was assumed to be constant.

\subsection{Convex multi-index transportation problems}

Our result on multi-index transportation problems and other applications
follow from a more general result about $n$-fold integer programming.
We begin this section with a few facts about this theory;
see \cite{DHOW,HKW,HOW,Onn} for more details.

Linear {\em $n$-fold integer programming} is the following problem in
dimension $nt$,
\begin{equation}\label{NIP}
\max\left\{wx\ :\ \A x=b\,,\ l\leq x\leq u\,,\ x\in\Z^{nt}\right\}\ ,
\end{equation}

where
\begin{equation}\label{NFold}
A^{(n)}\quad:=\quad
\left(
\begin{array}{cccc}
  A_1    & A_1    & \cdots & A_1    \\
  A_2    & 0      & \cdots & 0      \\
  0      & A_2    & \cdots & 0      \\
  \vdots & \vdots & \ddots & \vdots \\
  0      & 0      & \cdots & A_2    \\
\end{array}
\right)\quad
\end{equation}
is an $(r+ns)\times nt$ matrix which is the {\em $n$-fold product}
of a fixed {\em $(r,s)\times t$ bimatrix} $A$, that is, of a matrix $A$
consisting of two blocks $A_1$, $A_2$, with $A_1$ its $r\times t$
submatrix consisting of the first $r$ rows and $A_2$ its $s\times t$
submatrix consisting of the last $s$ rows.

Let $A$ be a fixed integer $(r,s)\times t$ bimatrix.
For any $n$ we write each vector $x\in\Z^{nt}$ as a tuple
$ x=(x^1,\dots, x^n)$ of $n$ {\em bricks} $x^i\in\Z^t$.
It has been shown in \cite{AT}, \cite{SS}, and \cite{HS}, in increasing
generality, that for every bimatrix $A$, the number of nonzero bricks
appearing in any element in the Graver basis $\G(\A)$
for any $n$ is bounded by a constant independent of $n$.
So we can make the following definition.
\bd{GraverComplexity}
The {\em Graver complexity} of an integer bimatrix $A$ is
defined to be the largest number $g(A)$ of nonzero bricks
$g^i$ in any element $ g\in\G\left(\A\right)$ for any $n$.
\ed

We then have the following lemma concerning Graver bases of $n$-fold products.
\bl{n-fold_edge}
For any integer bimatrix $A$ and all $n$ we have $\|G(\A)\|_1\leq g(A)\|G(A_2)\|_1$.
\el
\bpr
Consider any Graver basis element $g\in\G\left(\A\right)$ for any $n$.
Then $\A g= 0$ and hence $\sum_{i=1}^n A_1 g^i= 0$ and
$A_2 g^i= 0$ for all $i$. Therefore (see \cite[Chapter 3]{Onn}), for each $i$,
$g^i$ can be written as the sum $ g^i=\sum_{j=1}^{k_i} h^{i,j}$
of some elements $h^{i,j}\in\G(A_2)$ all lying in the same orthant.
Let $m:=k_1+\cdots+k_n$ and let $ h$ be the vector
$$ h\ :=\ ( h^{1,1},\dots, h^{1,k_1},\dots,
h^{n,1},\dots, h^{n,k_n})\ \in\ \Z^{mt}\ .$$
Then $\sum_{i,j}A_1 h^{i,j}= 0$ and $A_2 h^{i,j}= 0$
for all $i,j$ and hence $A^{(m)} h= 0$. We claim that
moreover, $ h\in\G\left(A^{(m)}\right)$. Suppose indirectly
this is not the case. Then there is an $\bar h\in\G\left(A^{(m)}\right)$
with $\bar h\sqsubset h$. But then the vector $\bar g\in\Z^{nt}$
defined by $\bar g^i:=\sum_{j=1}^{k_i}\bar h^{i,j}$ for all $i$
satisfies $\bar g\sqsubset\bar g$ contradicting $ g\in\G\left(\A\right)$.
This proves the claim. Therefore, by Definition \ref{GraverComplexity}
the number of nonzero bricks $ h^{i,j}$ of
$h$ is at most $g(A)$. So
$$\|g\|_1 = \sum_{i=1}^n\|g^i\|_1 =
\sum_{i=1}^n\left\|\sum_{j=1}^{k_i} h^{i,j}\right\|_1 \leq
\sum_{i=1}^n\sum_{j=1}^{k_i}\| h^{i,j}\|_1 \leq
m\,\|G(A_2)\|_1 \leq g(A)\|G(A_2)\|_1 .$$
Since this holds for all $g\in\G(\A)$, we get 
$\|G(\A)\|_1\leq g(A)\|G(A_2)\|_1$ for all $n$.
\epr

We now have the following theorem on convex
$n$-fold integer programming.

\bt{N-fold}
For every fixed $d$, $p$ and $(r,s)\times t$ bimatrix $A$,
there exists a smallest constant $v(d,p;A)$ such that,
for every $n$, every $b\in\Z^{r+ns}$, every $l,u\in\Z^{nt}$, and
every $\{0,1,\dots, p\}$-valued $d\times(nt)$ matrix $W$,
the projection $\conv(WS)$ of the set
\begin{equation}\label{N-fold_equation}
S\ :=\ \left\{x\in\Z^{nt}\,,\ \A x=b\,,\ l\leq x\leq u\right\},
\end{equation}
has at most $v(d,p;A)$ vertices. Further, all vertices of $\conv(WS)$
can be enumerated, and for any convex $f:\R^d\rightarrow\R$
the convex $n$-fold program $\max\{f(Wx)\,:\,x\in S\}$ solved, by
solving $v(d,p;A)$ linear counterparts over $S$, and in polynomial time.
\et
\bpr
Since $A$ and in particular its second block $A_2$ are fixed,
it follows from Lemmas \ref{Graver_projections} and \ref{n-fold_edge}
that $S$ has edge complexity
$e(S)\leq\|G(\A)\|_1\leq g(A)\|G(A_2)\|_1$ bounded by a constant.
So it follows from Theorem \ref{main} that the number of vertices of
$\conv(WS)$ is bounded by a constant $v(d,p;A)$ independent
of $n,b,l,u,W$, and that the convex problem reduces to solving
$v(d,p;A)$ linear counterparts over $S$. Now each linear counterpart
is a linear $n$-fold integer programming problem and is solvable in
time polynomial in $n$ and $\l b,l,u\r$, see \cite{DHOW,HOW},
implying the theorem.
\epr

The theory of $n$-fold integer programming has a variety of application,
including to multi-commodity flows, privacy in databases, and more, see
for example \cite[Chapters 4--5]{Onn}. Here we mention only one, but generic,
consequence of Theorem \ref{N-fold}, to multi-index transportation programs;
indeed such integer programs are universal: the universality theorem
of \cite{DO} asserts that any integer program is polynomial time liftable
to some isomorphic $3\times m\times n$ multi-index transportation program.

\bc{Multi-index}
For every fixed $d,p,l,m$, there exists a constant $t(d,p;l,m)$ such that,
for every $n$, every integer line-sums $a_{j,k},b_{i,k},c_{i,j}$,
and every $\{0,1,\dots, p\}$-valued $d\times(l\times m\times n)$ matrix $W$,
the projection $\conv(WS)$ of the set of $l\times m\times n$ tables
\begin{equation}\label{tables}
S\ :=\ \left\{x\in\Z_+^{l\times m\times n}\ :\ \sum_i x_{i,j,k}=a_{j,k}
\,,\ \sum_j x_{i,j,k}=b_{i,k}\,,\ \sum_k x_{i,j,k}=c_{i,j}\right\}\ ,
\end{equation}
has at most $t(d,p;l,m)$ vertices regardless of $n$ and the line sums.
Moreover, all these vertices can be enumerated, and for any convex $f$
the problem $\max\{f(Wx):x\in S\}$ solved, by solving $t(d,p;l,m)$
linear counterparts over $S$, and in polynomial time.
\ec
\bpr
Let $A$ be the $(lm,l+m)\times lm$ bimatrix whose first block $A_1$
is the $lm\times lm$ identity matrix $I_{lm}$ and second block $A_2$
is the $(l+m)\times lm$ vertex-edge incidence matrix of the complete
bipartite graph $K_{l,m}$. Then it is not hard to verify that, suitably
arranging tables in vectors in $\Z^{l\times m\times n}\cong\Z^{lmn}$
and the line-sums $a_{j,k},b_{i,k},c_{i,j}$ in a right-hand side vector
$e\in\Z^{lm+n(l+m)}\cong\Z^{lm+ln+mn}$, the line-sum equations
are encoded in the system $\A x=e$, encoding the
multi-index transportation program in a suitable $n$-fold program.
The corollary now follows from Theorem \ref{N-fold}.
\epr

The determination of the constants $v(d,p;A)$ in Theorem \ref{N-fold}
for bimatrices $A$ in various applications, and in particular the
constants $t(d,p;l,m)$ in Corollary \ref{Multi-index}, is
extremely difficult. Let us briefly discuss the case of binary
criteria matrices $W$ and define $t(d;l,m):=t(d,1;l,m)$ to be
the largest possible number of vertices of $\conv(WS)$ for any
$n$, any $\{0,1\}$-valued $d\times(l\times m\times n)$ matrix $W$,
and any set $S$ of $l\times m\times n$ tables as in (\ref{tables}).
We will bound it from above in terms of the relevant Graver complexity.
Now, the relevant bimatrix (see proof of Corollary \ref{Multi-index})
is $A$ with first block $A_1=I_{lm}$ and second block $A_2$
the $(l+m)\times lm$ vertex-edge incidence matrix of the complete
bipartite graph $K_{l,m}$. Let $g(l,m):=g(A)$ be the Graver complexity
of this matrix. We note that already $g(l,m)$ is not known for
all values with $3\leq l<m$, but bounds are available, \cite{KT}.
We have the following upper bound.

\bp{binary_table_upper_bound}
The maximum number of vertices of the projection $\conv(WS)$
of any set of $l\times m\times n$ tables $S$ in (\ref{tables}) by any
binary $d\times(l\times m\times n)$ matrix $W$ satisfies
$$t(d;l,m)\leq2\sum_{k=0}^{d-1}
{{{1\over2}\left((2g(l,m)\min\{l,m\}+1)^d-3\right)}\choose k}
=O\left((2g(l,m)\min\{l,m\})^{d(d-1)}\right) .$$
\ep
\bpr
The second block $A_2$ of $A$ is totally unimodular so its Graver basis
$\G(A_2)$ consist precisely of its circuits, see proof of Theorem \ref{TUM};
and the circuits of $A_2$ are vectors supported on circuits of $K_{l,m}$ with
alternating values of $1$ and $-1$. So each $h\in\G(A_2)$ consists of
the same number, at most $\min\{l,m\}$, of $-1$ and $1$ values.

Now, the proofs of Theorem \ref{main} and Lemma \ref{Graver_projections}
imply that every edge of $\conv(WS)$ is parallel to $Wg$ for some
element $g\in\G(\A)$. Consider any $k=1,\dots,d$ and let $W_k$ be the
corresponding row of $W$. Write $g$ and $W_k$ as tuples
$g=(g^1,\dots,g^n)$ and $W_k=(W_k^1,\dots,W_k^n)$ of bricks with each
$g^i$ and $W_k^i$ in $\Z^{l\times m}$. Now, for each $i$, $g^i$ can be
written as a sum $g^i=\sum_{j=1}^{k_i} h^{i,j}$ of elements $h^{i,j}\in\G(A_2)$,
all in the same orthant, with $\sum_{i=1}^nk_i\leq g(l,m)$
(see proof of Lemma \ref{n-fold_edge}). Since $W$ is binary we get
$$|W_kg|=\left|\sum_{i=1}^n\sum_{j=1}^{k_i} W_k^ih^{i,j}\right|\leq
\sum_{i=1}^n\sum_{j=1}^{k_i}\left|W_k^ih^{i,j}\right|\leq
\sum_{i=1}^n\sum_{j=1}^{k_i}\min\{l,m\}\leq g(l,m)\min\{l,m\}\ .$$

So $\conv(WS)$ is refined by a zonotope with
${1\over2}\left(\left(2g(l,m)\min\{l,m\}+1\right)^d-1\right)$ generators,
consisting of one element of each antipodal pair of nonzero vectors in
$$\left\{0,\pm 1,\pm 2, \dots,\pm g(l,m)\min\{l,m\}\right\}^d\ .$$
Therefore the claimed upper bound now follows from the bound in Lemma \ref{Zonotope}.
\epr

Of particular interest is the case of $3\times 3\times n$ tables, that is,
$l=m=3$ and arbitrary $n$, which is the smallest where the line-sum equations
defining $S$ in (\ref{tables}) are already not totally unimodular.
Here the Graver complexity is known to be $g(3,3)=9$. For
bicriteria optimization, that is, $d=2$, the bound on the largest possible
number of vertices of $\conv(WS)$ for any $n$, any $\{0,1\}$-valued
$2\times(3\times 3\times n)$ matrix $W$, and any set $S$ of  $3\times 3\times n$
tables with any line-sums as in (\ref{tables}), is $t(2;3,3)\leq 3024$.
It would be interesting to determine the value of $t(2;3,3)$ exactly.

\section{Open problems}

This article raises many open problems, in particular of
combinatorial-geometric flavor concerning the determination of
the various bounds on the number of vertices of projections
$\conv(WS)$ of various sets $S$ under small valued and binary $W$.

First, as discussed in Section 4, we know very little on the values $m(d,p)$.
While the relation $u(d)\leq m(d)\leq z(d)$ on the
maximum numbers $m(d)=m(d,1)$ and $u(d)$ of vertices of the projection of any
matroid and any uniform matroid $S$, respectively, under any binary $W$,
and the number $z(d)$ of vertices of the zonotope $Z(d)=\zone\{-1,0,1\}^d$,
holds with equality for $d=1$ trivially and for
$d=2$ with $u(2)=m(2)=z(2)=8$ by Theorem \ref{Bounds}, we know nothing
about higher $d$. The determination of $z(d)$ should certainly be
easier and in particular could be obtained by brute force computation
for small values of $d$. It would be particularly interesting
to determine $u(d)$ for all $d$, while $m(d)$ may be still way out of reach.

Turning to Section 5.1, an interesting open question related to integer
programs with totally unimodular defining systems concerns projections
of the so-called {\em assignment polytope} (often also termed
{\em Birkhoff polytope} or {\em bistochastic polytope}) which is the convex
hull $\conv(S)$ of the set $S$ of $n\times n$ permutation matrices given by
$$S\ =\ \{X\in\Z_+^{n\times n}\,:\,\ \sum_i X_{i,j}=1\,,\ \sum_j X_{i,j}=1\}\ .$$
Let $a_d(n)$ be the maximum number of vertices of
the projection $\conv(WS)$ of the assignment polytope by any
$\{0,1\}$-valued $d\times(n\times n)$ matrix $W$ into $\R^d$.
The assignment polytope has an exponential number
${1\over2}\sum_{k=2}^n{n\choose k}^2k!(k-1)!\geq{1\over n}{n!\choose2}$
of edge-directions, but by Theorem \ref{TUM}, for any fixed $d$,
we have that $a_d(n)$ is bounded by a polynomial in $n$. It is
interesting to determine $a_d(n)$ and in particular $a_2(n)$.

Finally, as discussed in Section 5.2, we know almost nothing even about
the maximum number $t(2;l,m)$ of vertices of the projection $\conv(WS)$
into the plane of any set of $l\times m\times n$ tables of any length
$n$ and any line-sums under the projection by any binary
$2\times(l\times m\times n)$ matrix $W$, except for the upper bound
$(2g(l,m)\min\{l,m\}+1)^2-1$ from
Proposition \ref{binary_table_upper_bound} in terms of the Graver
complexity $g(l,m)$. The determination of $g(l,m)$ for all $3\leq l<m$
is also open and challenging, with the smallest unknown value satisfying
$g(3,4)\geq 27$. Last, the value $t(2;3,3)$ for binary projections to
the plane of $3\times 3\times n$ tables is particularly intriguing.

\section*{Acknowledgments}
Shmuel Onn is supported in part by a grant from the Israel Science Foundation.
Michal Rozenblit is supported by a scholarship from the Technion Graduate School.

\end{document}